\documentclass[12pt]{amsart}

\usepackage{amsmath,amsfonts,amscd,amsthm,amssymb,latexsym}

\font\teneufm=eufm10 \font\seveneufm=eufm7 \font\fiveeufm=eufm5
\newfam\frakturfam
\textfont\frakturfam=\teneufm \scriptfont\frakturfam=\seveneufm
\scriptscriptfont\frakturfam=\fiveeufm

\newtheorem{pr}{Proposition}

\newtheorem{lm}{Lemma}
\newtheorem{theor}{Theorem}
\newtheorem{co}{Corollary}

\def\bee{\begin{eqnarray}}
\def\bes{\begin{eqnarray*}}
\def\eee{\end{eqnarray}}
\def\ees{\end{eqnarray*}}

\def\Proof{{\sl Proof.}\ }

\title{On the solvability of graded Novikov algebras}

\begin{document}
\date{}
\maketitle

\begin{center}

{\bf Ualbai Umirbaev}\footnote{Department of Mathematics,
 Wayne State University,
Detroit, MI 48202, USA; Department of Mathematics, 
Al-Farabi Kazakh National University, Almaty, 050040, Kazakhstan; 
and Institute of Mathematics and Mathematical Modeling, Almaty, 050010, Kazakhstan,
e-mail: {\em umirbaev@wayne.edu}}
and
{\bf Viktor Zhelyabin}\footnote{Institute of Mathematics of the SB of RAS, Novosibirsk, 630090, Russia,
e-mail: {\em vicnic@math.nsc.ru}}

\end{center}

\begin{abstract} We show that the right ideal of a Novikov algebra generated by  the square  of a right nilpotent subalgebra is nilpotent. We also prove that a $G$-graded Novikov algebra $N$ over a field $K$ with solvable $0$-component $N_0$ is solvable, where $G$ is a finite additive abelean group and the characteristic of $K$ does not divide the
order of the group $G$. We also show that any Novikov algebra $N$ with a finite solvable group of automorphisms $G$ is solvable if the algebra of invariants $N^G$ is solvable. 
\end{abstract}


\noindent


\section{Introduction}

\hspace*{\parindent}

A nonassociative algebra  $N$ over a field  $K$ is called a {\em Novikov} algebra \cite{Osborn1992} if
it satisfies the following  identities:
\bee\label{f1}
(x,y,z)=(y,x,z)\, \text{ (left symmetry)},
\eee
\bee\label{f2}
(xy)z=(xz)y\, \text{ (right commutativity)},
\eee
 where  $(x,y,z)=(xy)z-x(yz)$ is the associator of elements $x,y,z$.

The defining identities of a Novikov algebra first appeared in the study of Hamiltonian operators in the formal calculus of variations by I.M. Gelfand and
I.Ya. Dorfman \cite{gel_dor}. These identities played a crucial role in the classification of linear Poisson brackets of hydrodynamical type by A.A. Balinskii and  S.P. Novikov \cite{bal_nov}. 

In 1987 E.I. Zelmanov \cite{Zel} proved that all finite-dimensional simple Novikov algebras over a field $K$ of characteristic $0$ are one-dimensional.  V.T. Filippov \cite{Fil89} constructed a wide class of simple Novikov algebras of characteristic $p\geq 0$. 
 J.M. Osborn \cite{Osborn1992,Osborn1992CommAl,Osborn1994} and X. Xu \cite{Xu96,Xu01} continued the study of simple finite dimensional algebras over fields of positive characteristic and simple infinite dimensional algebras over fields of characteristic zero.  A complete classification of finite dimensional simple Novikov algebras over algebraically closed fields of characteristic $p>2$ is given in \cite{Xu96}.

E.I. Zelmanov  also proved that if $N$ is a finite
dimensional right nilpotent Novikov algebra then $N^2$ is nilpotent \cite{Zel}. In 2001 V.T. Filippov \cite{Fil01} proved that any left-nil Novikov algebra of bounded index over a field of characteristic zero is nilpotent. A.S. Dzhumadildaev and K.M. Tulenbaev \cite{DzT} proved that any right-nil Novikov algebra of bounded index $n$ is right nilpotent if the characteristic $p$ of the field $K$ is $0$ or $p>n$. In 2020 I. Shestakov and Z. Zhang proved \cite{ShZh} that for any
Novikov algebra $N$ over a field the following conditions are
equivalent:

$(i)$ $N$ is solvable;

$(ii)$ $N^2$ is nilpotent;

$(iii)$ $N$ is right nilpotent.

 The Freiheitssatz for Novikov algebras over fields of characteristic $0$ was proven by L. Makar-Limanov and U. Umirbaev \cite{MLU11N}. 
 L.A. Bokut, Y. Chen, and Z. Zhang \cite{BCZ} proved that every Novikov algebra is a subalgebra of a Novikov algebra obtained from some differential algebra by Gelfand-Dorfman construction \cite{gel_dor}.

This paper is devoted to the study of solvable, nilpotent, and right nilpotent Novikov algebras and graded Novikov algebras. 
Notice that  an  algebra $A$ over a field containing all $n$th roots of unity admits an automorphism
of order $n$ if and only if $A$ admits a $\mathbb{Z}_n$-grading. For this reason the study of graded algebras is related to the study of actions of finite groups. First we recall some definitions and classical results. 

 Let $R$ be an algebra over
a field $K$. For any automorphism $\phi$ of $R$ the set of fixed
elements 
\bes
R^{\phi} = \{x \in R|\phi(x) = x\}
\ees
 is a subalgebra of $R$
and is called the subalgebra of {\em invariants} of $\phi$. An
automorphism $\phi$ is called {\em regular} if $R^{\phi}=0$. For any group
$G$ of automorphisms
 of $R$ the subalgebra of invariants 
\bes
R^G = \{x \in  R|\phi(x) = x \text{ for
all }\phi \in G\}
\ees
 is defined similarly.

In 1957 G. Higman \cite{Higman} published a classical result on Lie
algebras which says that if a Lie algebra $L$ has a regular
automorphism  $\phi$ of prime order $p$, then $L$ is nilpotent. It
was also shown that the index of nilpotency $h(p)$ of $L$ depends
only on $p$. An explicit estimation of the function $h(p)$ was found
by A.I. Kostrikin and V.A. Kreknin \cite{Kreknin_Kostrikin} in 1963.
A little later, V.A. Kreknin proved \cite{Kreknin} that a finite
dimensional  Lie algebra  with a regular automorphism of an
arbitrary finite order is solvable.    In 2005 N. Yu. Makarenko
\cite{Makarenko} proved that if a Lie algebra $L$ admits an
automorphism of prime order $p$ with a finite-dimensional fixed
subalgebra of dimension $t$, then $L$ has a nilpotent ideal of
finite codimension with the index of nilpotency bounded in terms of
$p$ and the codimension bounded in terms of $t$ and $p$.

In 1973 G. Bergman and I. Isaacs \cite{BergIsac} published a
classical result on the actions of finite groups on associative
algebras. Let $G$ be a finite group of automorphisms of an
associative algebra $R$  and suppose that $R$ has no $|G|$-torsion.
If the subalgebra of invariants $R^G$ is nilpotent then the
Bergman-Isaacs  Theorem \cite{BergIsac} states that $R$ is also
nilpotent.  Since then a very large number of papers have been
devoted to the study of automorphisms of associative rings. The
central problem of these studies was to identify the properties of
rings that can be transformed from the ring of invariants to the
whole ring. In 1974 V. K. Kharchenko \cite{Kharchenko} proved if
$R^G$ is a PI-ring then $R$ is a PI-ring under the conditions of the
Bergman-Isaacs  Theorem.

The Bergman-Isaacs Theorem was partially generalized by W.S.
Martindale and S. Montgomery \cite{Martindale_Montgomery} in 1977 to
the case of a finite group of {\em Jordan} automorphisms, that is a
finite group of automorphisms of the adjoint Jordan algebra
$R^{(+)}$.

 An analogue  of
Kharchenko's result for Jordan algebras was proven by A. P. Semenov
\cite{Semenov} in 1991. In particular, A. P. Semenov proved that if
$J^G$ is  a solvable
 algebra  over a field of characteristic zero, then so is the Jordan algebra
 $J$. His proof uses a deep result by E.I. Zel'manov \cite{Zel'manov} which
 says that every Jordan nil-algebra of bounded index over a field of characteristic zero is solvable. If a Jordan algebra
$J$ over a field of characteristic not equal to $2,3$ admits an
automorphism $\phi$ of order $2$  with solvable $J^\phi$,
then $J$ is solvable \cite{Zhelyabin}.

In the case of alternative algebras one cannot expect that
nilpotency of the invariant subalgebra implies the nilpotency of the
whole algebra. There is an example (see \cite{Dorofeev, ZSSS}) of a
solvable non-nilpotent alternative  algebra with an automorphism of
order two such that its  subalgebra of invariants is nilpotent. A
combination of Semenov's result \cite{Semenov} and Zhevlakov's
theorem \cite{Zhevlakov} gives that, for an alternative algebra $A$ over a
field of characteristic zero, the
solvability of the algebra of invariants $A^G$ for a finite group
$G$ implies the solvability of $A$. It is also known
\cite{Smirnov} that if $A$ is an alternative algebra over a field
of characteristic not equal to $2$ with an automorphism $\phi$ of
order two, then the solvability of the algebra of invariants
$A^\phi$ implies the solvability of $A$.   In \cite{Goncharov} M.
Goncharov proved that an alternative $\mathbb{Z}_3$- graded
algebra $A=A_0\oplus A_1\oplus A_2$ over a field of characteristic
not equal to $2,3,5$ is solvable if $A_0$ is solvable.

It was shown in \cite{ZhU} for every $n$ of the form $n=2^k3^l$ that 
a $\mathbb{Z}_n$- graded Novikov  
\bes
N=N_0\oplus\ldots\oplus
N_{n-1}
\ees
 over a field of characteristic not equal to $2, 3$ is
solvable if $N_0$ is solvable.

In this paper we first prove that if $L$ is a right nilpotent subalgebra of a Novikov algebra $N$ then the right ideal of $N$ generated by $L^2$ is right nilpotent (Theorem \ref{t1}). This result gives a deeper explanation of the results on the nilpotency of $N^2$ mentioned above. The main result of the paper (Theorem \ref{t2}) says that if $N$ is a $G$-graded Novikov algebra with solvable $0$-component $N_0$, where $G$ is a finite additive abelian group, then $N$ is solvable. This result allows us to prove (Theorem \ref{t3}) that if $N$ is a Novikov algebra with solvable algebra of invariants $N^G$, where $G$ is a finite solvable group of automorphisms of $N$, then $N$ is solvable. Theorems \ref{t2} and \ref{t3} are formulated for fields of characteristic $0$ or positive characteristic $p$ that does not divide $|G|$. 
Notice that the solvability and the right nilpotency of Novikov algebras are equivalent by the result of I. Shestakov and Z. Zhang mentioned above.

The paper is organized as follows. In Section 2 we prove some identities and Theorem  \ref{t1}.  Sections 3--5 are devoted to the study of the structure of $\mathbb{Z}_n$-graded Novikov algebras. Theorems  \ref{t2} and  \ref{t3} are formulated and proven in Section 6.

\section{Right nilpotent subalgebras}

\hspace*{\parindent}

 The identities (\ref{f1}) and (\ref{f2}) easily imply the identities 
\bee\label{f3}
(xy,z,t)=(x,z,t)y
\eee
and
\bee\label{f4}
(x,yz,t)=(x,y,t)z.
\eee

Let $A$ be an arbitrary algebra. The powers of $A$ are defined inductively by $A^1=A$ and
\bes
A^m=\sum_{i=1}^{m-1}A^{i}A^{m-i}
\ees
 for all positive integers $m\geq 2$. The algebra $A$ is called {\em nilpotent}
if  $A^{m}=0$ for some positive integer  $m$.

The right powers of $A$ are defined inductively by $A^{[1]}=A$ and $A^{[m+1]}=A^{[m]}A$ for all integers $m\geq 1$. The algebra $A$ is called {\em right nilpotent} if there exists a positive integer $m$ such that $N^{[m]}=0$. In general, the right nilpotency of an algebra does not imply its nilpotency. This is also true in the case of Novikov algebras.

{\em Example 1}. \cite{Zel} Let  $N=Fa+Fb$ be a vector space of dimension 2. The product on $N$ is defined as
 $$ab=b,a^2=b^2=ba=0.$$
 It is easy to check that $N$ is a right nilpotent Novikov algebra, but not nilpotent.

The derived powers of $A$ are defined by $A^{(0)}=A$,
$A^{(1)}=A^2$, and $A^{(m)}=A^{(m-1)}A^{(m-1)}$ for all positive
integers $m\geq 2$. The algebra $A$ is called {\em solvable} if
$A^{(m)}=0$ for some positive integer  $m$. Every right nilpotent
algebra is solvable, and, in general, the converse is not true. But
 every solvable Novikov algebra is right
nilpotent \cite{ShZh}.

  It is well known that if $I$
  and $J$ are ideals of a Novikov algebra $N$, then $IJ$ is also an ideal of $N$.
Consequently, if $N$ is a Novikov algebra then $N^m$, $N^{[m]}$, and $N^{(m)}$ are ideals of $N$.
If $S$ is a subset of a Novikov algebra $N$,
then denote by $\langle S\rangle$ the right ideal of $N$ generated by $S$. Notice that if $I$ is a right
ideal of $N$, then $IS$ is a right ideal of $N$ for any subset $S\subseteq N$ by (\ref{f2}).

In any algebra we denote by $x_1x_2\ldots x_k$ the right normed product $(\ldots(x_1x_2)\ldots )x_k$ of elements
 $x_1,x_2,\ldots,x_k$. For any $x,y$ denote by 
 $x\circ y=xy+yx$  
 the Jordan product.

\begin{lm}\label{l1}  Any Novikov algebra satisfies the following identities:
\bee\label{f5}
a(bx_1\ldots x_t)=abx_1x_2\ldots x_t-\sum_{i=1}^t(a,b,x_i)x_1\ldots x_{i-1}x_{i+1}\ldots x_t
\eee
for each positive integer $t\geq 1$,
\bee\label{f6}
(ax_1\ldots x_s)\circ(bx_{s+1}\ldots x_{t})=(a\circ b)x_1x_2\ldots x_t-\sum_{i=1}^k(a,b,x_i)x_1\ldots x_{i-1}x_{i+1}\ldots x_t
\eee
for all nonnegative integers $0\leq s<t$, and
\bee\label{f7}
(ax_1\ldots x_s)\circ(bx_{s+1}\ldots x_{t})=a\circ(bx_1\ldots x_t).
\eee
\end{lm}
\Proof We prove (\ref{f5}) by induction on $t$. If $t=1$, then (\ref{f5}) is true by the definition of the associator.
By (\ref{f4}), we have
\bes
a(bx_1\ldots x_t)=a(bx_1\ldots x_{t-1})x_t-(a,bx_1\ldots x_{t-1},x_t)\\
=a(bx_1\ldots x_{t-1})x_t-(a,b,x_t)x_1\ldots x_{t-1}.
\ees
Using this and the induction proposition, we get
\bes
a(bx_1\ldots x_t)
=(abx_1x_2\ldots x_{t-1}-\sum_{i=1}^{t-1}(a,b,x_i)x_1\ldots x_{i-1}x_{i+1}\ldots x_{t-1})x_t\\
-(a,b,x_t)x_1\ldots x_{t-1}
=abx_1x_2\ldots x_t-\sum_{i=1}^t(a,b,x_i)x_1\ldots x_{i-1}x_{i+1}\ldots x_t.
\ees
By (\ref{f2}), (\ref{f3}), and (\ref{f5}), we get
\bes
(ax_1\ldots x_s)(bx_{s+1}\ldots x_{t})=(ab)x_1x_2\ldots x_t-\sum_{i=s+1}^t(a,b,x_i)x_1\ldots x_{i-1}x_{i+1}\ldots x_t.
\ees
This implies (\ref{f6}). The identity (\ref{f7}) is a direct consequence of (\ref{f2}), (\ref{f5}), and (\ref{f6}). $\Box$

Let $N$ be a Novikov algebra and let $L$ be a subalgebra of $N$.
Set $L_0=N$ and $L_k=\langle L^{[k]}\rangle$ for each positive integer $k$.

Consider the descending sequence of right ideals
\bes
N=L_0\supseteq L_1\supseteq L_2\supseteq\ldots \supseteq L_k\supseteq \ldots
\ees
of the algebra $N$.

\begin{lm}\label{l2}
$L_sL_t\subseteq L_{s+t-1}$ for all positive integers $s,t$.
\end{lm}
\Proof We prove the lemma by induction on $t$. It is true for $t=1$ by the definition of $L_s$. Notice that
\bes
L_s=L_1\underbrace{L\ldots L}_{s-1}
\ees
for each $s\geq 1$ by (\ref{f2}).

Suppose that $t\geq 2$ and let $x\in L_s$ and $y=za_1\ldots a_{t-1}\in L_t$, where $z\in L_1$ and $a_1,\ldots,a_{t-1}\in L$.
By (\ref{f5}), we get
\bes
xy=xza_1\ldots a_{t-1}-\sum_{i=1}^{t-1}(x,z,a_i)a_1\ldots \widehat{a_i}\ldots a_{t-1}
\ees
where $\widehat{a_i}$ means that $a_i$ is absent. Notice that $xz\in L_s$ and 
\bes
xza_1\ldots a_{t-1}\in L_s\underbrace{L\ldots L}_{t-1}=L_{s+t-1}.
\ees
Moreover,  $(x,z,a_i)$ belongs to the right ideal generated by $(L^{[s]},L,a_i)$ by (\ref{f3}) and (\ref{f4}). Consequently, $(x,z,a_i)\in L_{s+1}$ and
$(x,z,a_i)a_1\ldots \widehat{a_i}\ldots a_{t-1}\in L_{s+t-1}$.
$\Box$

In general, $L_1$ is not an ideal of $L_0=N$.

{\em Example 2}. Let $K[x,y]$ be the polynomial algebra over $K$
in the variables $x, y$. Define a new product $\cdot$ on $K[x,y]$ by  $$f\cdot g=f\frac{\partial
g}{\partial x},\, f,g \in K[x,y].$$ Then $N=(K[x,y],\cdot)$ is a
Novikov algebra.

Let $L=Kx$. Then $L$ is a subalgebra of $N$ since $x\cdot x=x$.
Let $L_1=\langle L\rangle$. It is clear that $L_1\subseteq
xK[x,y]$. Hence, $$y\cdot x =y\frac{\partial x}{\partial
x}=y\not\in L_1.$$ Consequently,  $L_1$ is not an ideal of
$L_0=N$.

But for each $r\geq 2$ the right ideal $L_r$ is an ideal of $L_1$ by Lemma \ref{l2}.

\begin{co}\label{c1}
$L_2^n\subseteq L_{n+1}$ for all $n\geq 1$.
\end{co}
\Proof It is trivial for $n=1$ and true for $n=2$ by Lemma \ref{l2}. If $L_2^i\subseteq L_{2+i-1}$ and $L_2^j\subseteq L_{2+j-1}$, then
$L_2^iL_2^j\subseteq L_{i+1}L_{j+1}\subseteq L_{i+j+1}$. Leading an induction on $n$ we get
\bes
L_2^n=\sum_{i+j=n, i,j\geq 1}L_2^iL_2^j\subseteq L_{n+1}. \ \ \ \Box
\ees

\begin{theor}\label{t1} Let $L$ be a right nilpotent subalgebra of a Novikov algebra $N$ over a field $K$.
Then the right ideal $L_2=\langle L^2\rangle$ of $N$ generated by $L^2$ is nilpotent.
\end{theor}
\Proof Suppose that $L^{[n]}=0$ for some $n\geq 2$. Then $L_n=0$. By Corollary \ref{c1}, we have $L_2^{n-1}\subseteq L_{n}=0$.
This means that $L_2$ is nilpotent. $\Box$

\section{$\mathbb{Z}_n$-graded Novikov algebras}

\hspace*{\parindent}

Let $\mathbb{Z}_n=\mathbb{Z}/n\mathbb{Z}$ be the additive cyclic group of order $n$.
Let
\bee\label{f8}
N=N_0\oplus N _1\oplus N_2\oplus\ldots\oplus N_{n-1}, \ \ N_iN_j\subseteq
N_{i+j}, \ i,j\in \mathbb{Z}_n,
\eee
 be a $\mathbb{Z}_n$-graded Novikov algebra over $K$.

If $f\in N_i$ then we say that $f$ is a homogeneous element of degree $i$.
Notice that $i$ is an element of $\mathbb{Z}_n$. Sometimes we consider the subscripts $i$ of $N_i$ as integers satisfying the condition 
$0\leq i\leq n-1$.

Obviously, $A=N_0$ is a subalgebra of $N$. Recall that $A^{[r]}$ is the right $r$th power of $A$.

\begin{lm}\label{l3}  Let $i_1,i_2,\ldots,i_k\in \mathbb{Z}_n$ and $i_1+i_2+\ldots+i_k=0$. Then
\bes
 A^{[r]}N_{i_1}N_{i_2}\ldots N_{i_k}\subseteq A^{[r]}.
\ees
\end{lm}
\Proof By the definition of a $\mathbb{Z}_n$-graded algebra, we have
\bes
AN_{i_1}N_{i_2}\ldots N_{i_k}\subseteq A.
\ees
Using this and (\ref{f2}), we get
\bes
A^{[r]}N_{i_1}N_{i_2}\ldots N_{i_k}=AN_{i_1}N_{i_2}\ldots N_{i_k}\underbrace{A\ldots A}_{r-1}\subseteq A\underbrace{A\ldots A}_{r-1} =A^{[r]}. \ \ \Box
\ees

Set $A^{\{0\}}=N$ and for any integer $r\geq 1$ denote by $A^{\{r\}}=\langle A^{[r]}\rangle$ the right ideal of $N$ generated by $A^{[r]}$.
Obviously, $A^{\{r\}}$ is a $\mathbb{Z}_n$-graded algebra, i.e.,
\bes
A^{\{r\}}=A^{\{r\}}_0\oplus A^{\{r\}} _1\oplus A^{\{r\}}_2\oplus\ldots\oplus A^{\{r\}}_{n-1}.
\ees

\begin{co}\label{c2}  If $r\geq 1$ and $0\leq i\leq n-1$, then
\bes A^{\{r\}}_i=\sum_{i_1,i_2,\ldots,i_k}
A^{[r]}N_{i_1}N_{i_2}\ldots N_{i_k}, \ees where $0\leq
i_1,i_2,\ldots,i_k\leq n-1$, $i_1+i_2+\ldots+i_k\equiv i (\mathrm{mod} \ n)$ and
$i_1+i_2+\ldots+i_k<n$.

In particular, $A^{\{r\}}_0=A^{[r]}$.
\end{co}
Consider the descending sequence of right ideals
\bee\label{f9}
N=A^{\{0\}}\supseteq A^{\{1\}}\supseteq\ldots \supseteq A^{\{r\}}\supseteq \ldots
\eee
of the algebra $N$ and the quotient algebra
\bee\label{f10}
B=A^{\{1\}}/A^{\{2\}}=B_0\oplus B_1\oplus B_2\oplus\ldots\oplus B_{n-1}, \ \ B_iB_j\subseteq
B_{i+j}, \ i,j\in \mathbb{Z}_n.
\eee
Notice that $B$ is  a right $N$-module. We establish some properties of the algebra $B$.

\begin{lm}\label{l4}  Let $B$ be the Novikov algebra defined by (\ref{f10}). Then

$(i)$ $B_0=A/A^2$;

$(ii)$ $BB_0=0$;

$(iii)$ $x\circ y=xy+yx=0$ for any $x\in B_i$ and $y\in B_{n-i}$.

\end{lm}
\Proof The statement (i) is true since $A^{\{r\}}_0=A^{[r]}$ by Corollary \ref{c2}. The statement (ii) is a direct corollary of the inclusion
$A^{\{r\}}A\subseteq A^{\{r+1\}}$.

Let $x=ax_1x_2\ldots x_s\in A^{\{1\}}_i$ and
$y=bx_{s+1}x_{s+2}\ldots x_t\in A^{\{1\}}_{n-i}$, where $a,b\in A$
and $x_r\in N_{k_r}$ for all $1\leq r\leq t$. If $i=0$, then
$A_i=A_{n-i}=A$ and $xy,yx\in A^2$. Suppose that $i, n-i\neq 0$.
Then $\Sigma_{r=1}^s k_r=i\neq 0$, $\Sigma_{r=s+1}^t l_r=n-i \neq
0$, and $\Sigma_{r=1}^t k_r=0$. In particular, $t>s\geq 1$. By
(\ref{f7}), we have \bes x\circ y=(ax_1\ldots
x_s)\circ(bx_{s+1}\ldots x_{t})=a\circ(bx_1\ldots x_t). \ees The
condition $\Sigma_{r=1}^t k_r=0$ implies that $bx_1\ldots x_t\in A$.
Consequently, $x\circ y\in A^2$. This proves $(iii)$. $\Box$

\section{Right nilpotency modulo $A^{\{1\}}$}

\hspace*{\parindent}

In this section we show that if the $0$-component $A=N_0$ of a $\mathbb{Z}_n$-graded Novikov algebra $N$ of the form (\ref{f8}) is right nlpotent, then
$N^{[m]}\subseteq A^{\{1\}}$ for some positive integer $m$.

\begin{lm}\label{l5}  Let $N$ be an arbitrary Novikov algebra and let $V$ be a subspace of $N$. Then for any $r\geq 1$ we have
\bes
NV^{[r]}V \subseteq \langle V^{[r]}\rangle+NV^{[r+1]}.
\ees
\end{lm}
\Proof By (\ref{f1}), we get
\bes
(NV^{[r]})V\subseteq (N,V^{[r]},V)+NV^{[r+1]}\subseteq (V^{[r]},N,V)+NV^{[r+1]} \subseteq  \langle V^{[r]}\rangle+NV^{[r+1]}. \ \ \Box
\ees

\begin{co}\label{c3}  If $r\geq 1$, then
\bes
N\underbrace{V\ldots V}_{r+1}\subseteq \langle V^{[r]}\rangle+NV^{[r+1]}.
\ees
\end{co}
\Proof It is true for $r=1$ by Lemma \ref{l5}. If it is true for some $r\geq 1$, then we get
\bes
N\underbrace{V\ldots V}_{r+2} \subseteq \langle V^{[r]}\rangle V+N V^{[r+1]}V \subseteq \langle V^{[r+1]}\rangle+NV^{[r+2]}
\ees
by (\ref{f2}) and Lemma \ref{l5}. $\Box$

\begin{lm}\label{l6} Let $N$ be an arbitrary $\mathbb{Z}_n$-graded Novikov algebra $N$ from (\ref{f8}) and suppose that the $0$-component $A=N_0$
of $N$ is right nilpotent. Then there exists a positive integer $m$ such that $N^{[m]}\subseteq A^{\{1\}}$.
\end{lm}
\Proof Suppose that $A^{[r]}=0$ for some positive integer $r$. By Corollary \ref{c3},
\bes
N\underbrace{A\ldots A}_{r+1}\subseteq \langle A^{[r]}\rangle+NA^{[r+1]}=0.
\ees
Again, by Corollary  \ref{c3}, we get
\bes
N\underbrace{N_i\ldots N_i}_{n}\subseteq \langle N_i^{[n-1]}\rangle+NN_i^{[n]}.
\ees
Notice that $N_i^{[n]}\subseteq A$. Consequently,
\bes
N\underbrace{N_i\ldots N_i}_{n+1}\subseteq (\langle N_i^{[n-1]}\rangle+NA)N_i\subseteq \langle N_i^{[n]}\rangle+NN_iA\subseteq  A^{\{1\}}+NN_iA.
\ees
Using this, we can easily show by induction on $s$ that
\bes
N\underbrace{N_i\ldots N_i}_{sn+1}\subseteq  A^{\{1\}}+NN_i\underbrace{A\ldots A}_{s}.
\ees
Cosequently,
\bes
N\underbrace{N_i\ldots N_i}_{(r+1)n+1}\subseteq  A^{\{1\}}+NN_i\underbrace{A\ldots A}_{r+1} \subseteq  A^{\{1\}}
\ees
since $N\underbrace{A\ldots A}_{r+1}=0$.

Thus, every $N_i$ acts on $N$ nilpotently modulo $A^{\{1\}}$ from the right hand side. Moreover, by (\ref{f2}), this action is commutative.
This easily implies the existence of an integer $m$ such that $N^{[m]}\subseteq A^{\{1\}}$. $\Box$

\section{Right nilpotency of $B$}

\hspace*{\parindent}

In this section we prove that any $\mathbb{Z}_n$-graded Novikov algebra $B$ defined by (\ref{f10}) is right nilpotent if the characteristic of $K$ does not divide $n$. Suppose that $N$ is a $\mathbb{Z}_n$-graded Novikov algebra of the form (\ref{f8}) satisfying the conditions

$(a)$ $NA=0$ and

$(b)$ $x\circ y=xy+yx=0$ for any $x\in N_i$ and $y\in N_{n-i}$ and for any $i\in \mathbb{Z}_n$.

All statements in this section are formulated for the algebra $N$.

First we prove the following lemma.
\begin{lm}\label{l7} Let $x\in N_{n-i}$, $u\in N_i^{[k]}$, $i\in \mathbb{Z}_n$, and $k\geq 1$.
Then $xu=-kux$.
\end{lm}
\Proof We prove the statement of the lemma by induction on $k$. If $k=1$, then it is true by $(b)$. Suppose that $k>1$ and $u=vy$, where
$v\in N_i^{[k-1]}$ and $y\in N_i$. Using (\ref{f1}), (\ref{f2}), and the induction proposition, we get
\bes
xu=x(vy)=-(x,v,y)+(xv)y=-(v,x,y)-(k-1)(vx)y\\
=-(vx)y+v(xy) -(k-1)(vx)y= -k(vy)x+v(xy)=-kux+v(xy).
\ees
Notice that $xy\in N_{n-i}N_i\subseteq A$ and $v(xy)=0$ by the condition $(a)$. Consequently, $xu=-kux$. $\Box$

\begin{co}\label{c4}  If the characteristic of the field $K$ does not divide $n$, then $N_i^{[n]}N_{n-i}=0$  for any $i\in \mathbb{Z}_n$.
\end{co}
\Proof Notice that $N_i^{[n]}\subseteq A$ and $N_{n-i}N_i^{[n]}=0$ by the condition $(a)$. Then Lemma \ref{l7} gives that
$nN_i^{[n]}N_{n-i}=0$. If the characteristic of $K$ does not divide $n$, then this gives $N_i^{[n]}N_{n-i}=0$. $\Box$

\begin{lm}\label{l8} If the characteristic of the field $K$ does not divide $n$, then
\bes
N\underbrace{N_i\ldots N_i}_{2n}=0
\ees
 for any $i\in \mathbb{Z}_n$.
\end{lm}
\Proof Corollary  \ref{c3} and the condition $(a)$ give that
\bes
N\underbrace{N_i\ldots N_i}_{n}\subseteq \langle N_i^{[n-1]}\rangle+NN_i^{[n]}\subseteq \langle N_i^{[n-1]}\rangle
\ees
since $N_i^{[n]}\subseteq A$. Notice that $i(n-1)=-i=n-i$ in $\mathbb{Z}_n$. This means $N_i^{[n-1]}\subseteq N_{n-i}$.
Cosequently,
\bes
N\underbrace{N_i\ldots N_i}_{n}\subseteq  \langle N_{n-i}\rangle.
\ees
Using (\ref{f2}) and $(a)$, we get
\bes
N\underbrace{N_i\ldots N_i}_{n+1}\subseteq  \langle N_{n-i}\rangle N_i = \langle N_{n-i}N_i\rangle =\langle N_iN_{n-i}\rangle.
\ees
Then
\bes
N\underbrace{N_i\ldots N_i}_{2n}\subseteq =\langle N_iN_{n-i}\rangle\underbrace{N_i\ldots N_i}_{n-1}=\langle N_i^{[n]}N_{n-i}\rangle.
\ees
Corollary \ref{c4} implies the statement of the lemma. $\Box$

\begin{pr}\label{p1} Let $N$ be a $\mathbb{Z}_n$-graded Novikov algebra of the form (\ref{f8}) satisfying the conditions

$(a)$ $NA=0$ and

$(b)$ $x\circ y=xy+yx=0$ for any $x\in N_i$ and $y\in N_{n-i}$ and for any $i\in \mathbb{Z}_n$.

 If the characteristic of the field $K$ does not divide $n$, then $N$ is right nilpotent.
\end{pr}
\Proof By Lemma \ref{l8}, every $N_i$ acts nilpotently on the right $N$-module $N$. Moreover, this action is commutative by (\ref{f2}). Consequently, $N$ acts nilpotently on $N$. $\Box$

\section{Solvability and right nilpotency}

\hspace*{\parindent}

The solvability and the right nilpotency of Novikov algebras are equivalent \cite{ShZh}. In this section we use these notions as synonyms.

\begin{pr}\label{p2} Let $N$ be a $\mathbb{Z}_n$-graded Novikov algebra of the form (\ref{f8}) such that $A=N_0$ is solvable.
 If the characteristic of the field $K$ does not divide $n$, then $N$ is solvable.
\end{pr}
\Proof Consider the descending sequence of right ideals (\ref{f9}).
By Lemma \ref{l6} there exists a positive integer $m$ such that
$N^{[m]}\subseteq A^{\{1\}}$. The algebra $B$ from (\ref{f10})
satisfies all conditions of Proposition \ref{p1} by Lemma \ref{l4}.
By Proposition \ref{p1} there exists a positive integer $t$ such
that $B^{[t]}=0$. This means that $(A^{\{1\}})^{[t]}\subseteq
A^{\{2\}}$. By Theorem \ref{t1}, the algebra $A^{\{2\}}$ is
nilpotent. Consequently, $A^{\{1\}}$ and $N$ are both solvable.
$\Box$

Let $G$ be an additive abelian group. We say that
\bes
N=\bigoplus_{g\in G} N_g
\ees
 is a $G$-graded algebra if $N_gN_h\subseteq N_{g+h}$ for all $g,h\in G$.

\begin{theor}\label{t2} Let $G$ be a finite additive abelian group and let $N$ be a $G$-graded Novikov algebra with solvable $0$-component $N_0$.
If the characteristic of the field $K$ does not divide the order of the group $G$, then $N$ is solvable.
\end{theor}

\Proof We prove the statement of the theorem by induction on the order $|G|$ of $G$. If $G=\mathbb{Z}_{n}$, then $N$
is solvable by Proposition \ref{p2}. 

Every finite abelian group is a direct sum of cyclic subgroups. 
Suppose that $G=\mathbb{Z}_{n_1}\oplus
\mathbb{Z}_{n_2} \oplus\ldots\oplus \mathbb{Z}_{n_k}$, where $n_i> 1$ for all $i$ and $k\geq 2$. Then $G=\mathbb{Z}_{n_1}
\oplus G_1$, where $G_1=\mathbb{Z}_{n_2}\oplus\ldots\oplus \mathbb{Z}_{n_k}$. 
Denote by $\mathrm{pr}$ the  projection of $G$ onto the group
$\mathbb{Z}_{n_1}$. Set 
\bes
N_i'=\sum_{g\in G,pr(g)=i}N _g, 
\ees
where
$i=0,1,\ldots, n_1-1$. It is easy to show that 
\bes
N=N'_0\oplus\ldots
N'_{n_1-1}
\ees
and $N$ is a $\mathbb{Z}_{n_1}$-graded algebra.

It is also clear that 
\bes
N_0'=\sum_{g\in G,pr(g)=0}N _g
\ees
 is a $G_1$-graded algebra and the $0$-component of $N_0'$ is
$N_0$. Since $|G_1|<|G|$ it follows that $N_0'$ is solvable by the induction proposition. 
Now we can apply Proposition \ref{p2} to the $\mathbb{Z}_{n_1}$-graded algebra $N$. 
 Hence $N$
is solvable. $\Box$

The statement of the next lemma is well known. 
\begin{lm} \label{l9} Let $G$ be a group of automorphisms of an arbitrary algebra $A$ and let $H$ be a normal subgroup of $G$. Then $A^H$ is $G$-invariant, the quotient group $G/H$ acts on $A^H$ by automorphisms, and $(A^H)^{G/H}=A^G$.
\end{lm}
\Proof Let $a\in A^H$ and let $g\in G$. Then $ghg^{-1}\in H$ for any
$h\in H$. Therefore, $a^{ghg^{-1}}=a$ and 
$$(a^g)^h=a^{gh}=a^{ghg^{-1}g}=(a^{ghg^{-1}})^g=a^g.$$ Consequently, the algebra $A^H$ is $G$-invariant. 
Let $g\in G$ and let
$\overline{g}$ be the image of  $g$ in $G/H$. Then $\overline{g}$
defines an automorphism of the algebra  $A$ by the rule
$a^{\overline{g}}=a^g$. This action is well defined. Hence
 the quotient group $G/H$ acts on $A^H$. It is easy to check that $(A^H)^{G/H}=A^G$. $\Box$

\begin{co} \label{c5} Let $N$ be a  Novikov algebra and let $G$ be a
finite abelian group of automorphisms of $N$. If the algebra $N^G$
is solvable and the characteristic of the field $K$ does not divide
the order of the group $G$, then $N$ is solvable.
\end{co}
\Proof  We may assume that  $K$ is algebraically closed. We prove the statement of the corollary by induction on the order $|G|$ of  $G$. 
If $G$ is a
simple group, then  $G\cong \mathbb{Z}_p$, where $p$ is a prime
number. Let $\phi$ be a generating element  of the group $G$. Then
$\phi^p=e$, where $e$ is the identity element of $G$. Let $\epsilon$ be a primitive $p$th root of unity and 
 let
$N_i=\ker(\phi-\epsilon^i) $ for all $0\leq i\leq p-1$. The indexes $i$ may be considered as elements of $\mathbb{Z}_p$ since $\epsilon^p=1$. 
Obviously, 
\bes
N=N_0\oplus\ldots\oplus N_{p-1}
\ees
and it is easy to check that $N_iN_j\subseteq N_{i+j}$ for all $i,j\in \mathbb{Z}_p$, i.e., $N$ is a $\mathbb{Z}_p$-graded algebra. 
 Moreover, $N_0=N^G$. By Proposition \ref{p2}, $N$ is solvable.

Let $H$ be a proper subgroup of $G$.
Then, by Lemma \ref{l9},   the quotient group $G/H$ acts on  $N^H$ by automorphisms and  $(N^H)^{G/H}=N^G$. 
We get that $N^H$ is solvable by the induction proposition since $|G/H|<|G|$. Now we can apply the induction proposition to the group $H$ and get that $N$ is solvable. $\Box$

\begin{theor} \label{t3} Let $N$ be a Novikov algebra and let $G$ be a
finite solvable  group of automorphisms of $N$. If the algebra $N^G$
is solvable and the characteristic of the field $K$ does not divide
the order of the group $G$, then $N$ is solvable.
\end{theor}
\Proof We prove the statement of the theorem by induction on $|G|$. 
The case of abelian groups is considered in Corollary \ref{c5}. Suppose that $G$ is not abelian. 
Then the commutator subgroup $G'$ of the solvable finite group $G$ is a proper normal subgroup.

By Lemma \ref{l9}, $(N^{G'})^{G/G'}=N^G$. Then the algebra $N^{G'}$ is solvable by the induction proposition since $|G/G'|<|G|$. 
Applying the induction proposition to $G'$, we get that $N$ is solvable. $\Box$

\end{document}